\documentclass[twocolumn,letterpaper,10pt,conference]{ieeeconf}

\IEEEoverridecommandlockouts                              
\usepackage{graphicx} % Required for inserting images
\usepackage{quiver}
\usepackage{verbatim}
\usepackage{todonotes}
\usepackage{xcolor}

\usepackage{amsmath, amsfonts, amssymb, amsthm, mathtools}
\let\labelindent\relax
\usepackage{enumitem}
\usepackage[ruled]{algorithm2e}
\usepackage{thmtools}
\usepackage{fancyhdr}
\usepackage{hyperref}
\usepackage{tikz}
\usepackage{comment}
\usepackage[normalem]{ulem}
\usepackage{booktabs}
\usepackage{subcaption}      % For subfigures and subcaptions
\usepackage[font=footnotesize]{caption}         % Enhances caption control (optional but recommended)

\usepackage{hyperref}
\usepackage[noabbrev,capitalize]{cleveref}
\crefformat{equation}{Eq.~#1}
\crefformat{figure}{Fig.~#1}
\usepackage{cite}

\newcommand{\maps}{\colon}
\newcommand{\R}{\mathbb{R}}
\newcommand{\Rbar}{\overline{\R}}
\newcommand{\argmin}{\mathrm{argmin}}

\newcommand{\id}{\mathrm{id}}
\newcommand{\define}[1]{\textbf{#1}}

\newcommand{\FF}{\mathcal{F}}

\definecolor{gatas-blue}{rgb}{0.3803921568627451,0.5333333333333333,0.6980392156862745}
\definecolor{devils-blue}{rgb}{1, 33, 105}

\usepackage{stmaryrd}
\DeclareMathOperator{\face}{\trianglelefteqslant}

\DeclareMathOperator{\image}{im}
\newcommand{\minimize}{\operatorname*{minimize}}

\renewcommand{\leq}{\leqslant}

\newtheorem{theorem}{Theorem}

\theoremstyle{definition}
\newtheorem{definition}{Definition}
\newtheorem*{interpretation}{Interpretation}
\newtheorem{example}{Example}
\newtheorem{examples}[example]{Examples}

\newtheorem{remark}{Remark}
\newtheorem{assumption}{Assumption}

\title{Distributed Multi-agent Coordination over Cellular Sheaves}
\author{Tyler Hanks, Hans Riess, Samuel Cohen, Trevor Gross, Matthew Hale, James Fairbanks
\thanks{$^*$Tyler Hanks, Samuel Cohen, and Trevor Gross are with the Department of Computer and Information Science and Engineering and James Fairbanks is with the Department of Mechanical and Aerospace Engineering, University of Florida, emails: \texttt{\{t.hanks,samuel.cohen,trevorgross,fairbanksj\}@ufl.edu}.
Hans Riess and Matthew Hale are with the Deparment of Electrical and Computer Engineering, Georgia Institute of Technology, emails: \texttt{\{riess, mhale30\}@gatech.edu}}
\thanks{Hanks was supported by the National Science Foundation Graduate Research Fellowship Program under Grant No. DGE-1842473. Any opinions, findings, and conclusions or recommendations expressed in this material are those of the author(s) and do not necessarily reflect the views of the NSF. Riess and Hale were supported by AFOSR under grants FA9550-23-1-0120 and FA9550-19-1-0169, ONR under grants N00014-22-1-2435 and N00014-21-1-2495. All authors were supported by DARPA under grant HR00112220038.}}

\date{\today}

\begin{document}

\maketitle

\begin{abstract}
    Techniques for coordination of multi-agent systems are vast and varied, often utilizing purpose-built solvers or controllers with tight coupling to the types of systems involved or the coordination goal. In this paper, we introduce a general unified framework for heterogeneous multi-agent coordination using the language of cellular sheaves and nonlinear sheaf Laplacians, which are generalizations of graphs and graph Laplacians. Specifically, we introduce the concept of a nonlinear homological program encompassing a choice of cellular sheaf on an undirected graph, nonlinear edge potential functions, and constrained convex node objectives, which constitutes a standard form for a wide class of coordination problems. We use the alternating direction method of multipliers to derive a distributed optimization algorithm for solving these nonlinear homological programs. To demonstrate the applicability of this framework, we show how heterogeneous coordination goals including combinations of consensus, formation, and flocking can be formulated as nonlinear homological programs and provide numerical simulations showing the efficacy of our distributed solution algorithm.
\end{abstract}

\section{Introduction}

Optimal control and coordination of multi-agent systems is of paramount importance in modern applications ranging from multi-robot systems to cyber-physical systems. In anticipation of next-generation systems composed of highly autonomous heterogeneous agents with multiple and evolving individual and group objectives, our goal in this paper is to present a unified framework for modeling a wide array of multi-agent coordination problems, along with a general-purpose solver for these problems. In this framework, system designers specify requirements about the behavior of a multi-agent system, compile these requirements into a \emph{homological program}, and synthesize a decentralized control policy with a universal solver. We make steps towards this goal.

The main abstractions we utilize in constructing this framework are cellular sheaves and nonlinear sheaf Laplacians. Our proposed framework uses cellular sheaves with nonlinear edge potentials to model various types of multi-agent interactions and coordination goals, including but not limited to, consensus, translation invariant formation, and flocking. The generality of this framework enables heterogeneity not only in the systems being controlled, but also in the communication patterns and coordination goals between different agents. To design distributed optimal controllers for such multi-agent systems, we use the alternating direction method of multipliers (ADMM) to derive a distributed solution algorithm for optimization problems defined over cellular sheaves. ADMM is a widely used technique for distributed optimization due to its attractive convergence properties provided only mild assumptions on the problem. The resultant algorithm allows agents to choose control inputs that minimize their own control objectives while driving the global system towards the coordination goal, all while utilizing only local computations and communication with neighboring agents.

%---------------------------
\paragraph*{Related Work}
%---------------------------
Several efforts have emerged to systematically study controller design patterns, including layered control architectures \cite{matni_quantitative_2024}, co-design, \cite{zardini_co-design_2021}, and information structures \cite{mahajan_information_2012}. Through our reformulation of well-known multi-agent coordination problems as homological programs, we similarly illuminate fundamental design principles of networked multi-agent systems.

The graph-theoretic approach to modeling multi-agent systems has lead to a corpus of algorithms addressing a diversity of coordination problems, including consensus \cite{olfati-saber_consensus_2004,olfati-saber_consensus_2007}, dissensus \cite{franci_breaking_2023}, formation control \cite{egerstedt_formation_2001,fax_information_2004}, flocking \cite{tanner_stable_2003,olfati-saber_flocking_2006}, and many other tasks  \cite{mesbahi_graph_2010}. In this literature, distributed optimization \cite{zheng_review_2022} is commonly integrated with cooperative control strategies in order to simultaneously pursue individual goals as well as shared common objectives \cite{raffard_distributed_2004}.

Cellular sheaves were introduced in stratification theory \cite{shepard_cellular_1985} and rediscovered in applied topology \cite{curry_sheaves_2014}. The proposed marriage of methods from spectral graph theory and cellular sheaves \cite{hansen_toward_2019} led to the observation that sheaf Laplacians induce a dynamical system generalizing (multidimensional) distributed mean consensus protocols where contemporaneous generalizations of consensus such as matrix-weighted \cite{trinh_matrix-weighted_2018} and matrix-scaled consensus \cite{trinh_matrix-scaled_2022} are straightforward instantiations of the cellular sheaf construction. Subsequent sheaf-theoretic developments \cite{ghrist_cellular_2022} yield a similar generalization of max consensus \cite{trinh_matrix-scaled_2022}.

Control of cellular sheaves was initiated for non-physical systems in which linear and nonlinear sheaf Laplacians are deployed to model opinion dynamics via discourse sheaves~\cite{hansen_opinion_2021}. Another sheaf-theoretic approach to optimal network control \cite{kearney_sheaf-theoretic_2020} bears similarities to ours, without the computational advantages of a distributed implementation. Sheaves have appeared elsewhere in the control systems literature, including modeling event-based systems \cite{zardini_compositional_2021}, information diffusion \cite{riess_diffusion_2022}, game theory \cite{riess_max-plus_2023}, as well as network theory \cite{robinson_understanding_2013}. Distributed optimization \cite{hanks2024compositional}, LQR \cite{she2023characterizing}, and model-predictive control \cite{hanks2023compositional} have also been studied through the lens of category theory, an area closely related to sheaf theory.

Homological programming was developed for mobile sensor networks~\cite{ghrist_positive_2017}, and extended to a generalized notion of distributed optimization~\cite{hansen_distributed_2019}, with applications in distributed model-predictive control and graph signal processing. Our work expands the scope of homological programming introduced in that work with the nonlinear sheaf Laplacian and constrained node objectives.

%---------------------------
\paragraph*{Outline}
%---------------------------

In \cref{sec:prelim}, we review cellular sheaves and sheaf Laplacians. Then, in \cref{sec:homological-programming}, we establish a general framework for representing multi-agent coordination problems. In \cref{sec:ADMM}, we propose a solution method. To encourage the adoption of this framework,  we present an in-depth case study in \cref{sec:examples}. In \cref{sec:simulations}, we present several numerical simulations. Finally, in \cref{sec:conclusion} we discuss future directions.

%%%%%%%%%%%%%%%%%%%%%%%%%%%%%%%%%%%%%%%%%%%
\section{Preliminaries}
\label{sec:prelim}
%%%%%%%%%%%%%%%%%%%%%%%%%%%%%%%%%%%%%%%%%%%

In this section we introduce cellular sheaves and sheaf Laplacians. Together with optimization techniques, these constructions will serve as a general framework for us to analyze multi-agent coordination problems.

Suppose $G = (V,E)$ is a undirected graph with node set $V = [N]$ and edge set $E \subseteq V \times V$, where $[N]$ denotes the set $\{1,2,\dots,N\}$. An edge between given nodes $i$ and $j$ is given by an unordered pair of concatenated indices $ij = ji \in E$. When particular endpoints of an edge are not specified, an edge can be denoted by $e \in E$. The set $N_i = \{j : ij \in E\}$ are the neighbors of $i$.

\begin{interpretation}
    In a multi-agent system (MAS), an undirected graph models bidirectional communication channels between agents. We assume these communication networks have unchanging topologies.
\end{interpretation}

Given a fixed graph $G$, we define a data structure for organizing assignments of vectors.

\begin{definition}[Cellular Sheaf]
Given $G = (V,E)$, a cellular sheaf $\mathcal{F}$ (valued in Euclidean spaces over a graph) is a data structure that assigns:
\begin{itemize}[leftmargin=*]
    \item A Euclidean space $\mathcal{F}(i)$ with the standard inner-product to every node $i \in V$ called a \define{stalk}.
    \item A Euclidean space $\mathcal{F}(ij)$ with the standard inner-product to every edge $ij \in E$ called an \define{edge stalk}.
    \item A linear transformation $\mathcal{F}(i) \xrightarrow{\mathcal{F}_{i \face ij}} \mathcal{F}(ij)$ for every $i \in V$, $j \in N_i$ called a \define{restriction map}.
\end{itemize}
\end{definition}

\begin{example}[Constant Sheaf]
    Given a graph $G$ and a particular vector space $\mathbb{R}^{k}$, the constant sheaf, denoted $\underline{\mathbb{R}}^k$, assigns the vector space $\underline{\mathbb{R}}^k(i) = \mathbb{R}^k$ to every $i \in V$, the vector space $\underline{\mathbb{R}}^k(ij) = \mathbb{R}^k$ to every edge $ij \in E$, and the identity map $\underline{\mathbb{R}}^k_{i \face ij} = \id_{\mathbb{R}^k}$ for every $i \in V$, $j \in N_i$.
\end{example}

We refer to a particular assignment $x_i \in \mathcal{F}(i)$ as a \define{local section}. This is the datum assigned to an individual node. Collecting all the data over nodes yields the space of $0$-cochains, $C^0(G; \mathcal{F}) = \bigoplus_{i \in V} \mathcal{F}(i)$, and data over edges is similarly gathered in the space of $1$-cochains, $C^1(G; \mathcal{F}) = \bigoplus_{e \in E} \mathcal{F}(e)$. $C^0(G;\mathcal{F})$ and $C^1(G;\mathcal{F})$ are endowed with the inner products $\langle \mathbf{x}, \mathbf{x}' \rangle_{C^0} = \sum_{i \in V} \langle x_i, x_i' \rangle_{\mathcal{F}(i)}$ and $\langle \mathbf{y}, \mathbf{y}' \rangle_{C^1} = \sum_{e \in E} \langle y_e, y_e' \rangle_{\mathcal{F}(e)}$.

\begin{interpretation}
    A stalk is interpreted as an agent's state space, and a local section is interpreted as a particular state (e.g.~position and velocity measurements).
    In a multi-agent systems, a $0$-cochain is, then, a global state, and $C^0(G;\mathcal{F})$ is the global state space. Note that dimensions of stalks are not required to be the same for every agent, hence, the state spaces of agents can vary from agent to agent (e.g.~a team of both ground and ariel vehicles).
    % In certain cases, $1$-cochains can be seen as observed outputs when internal states are transformed by restriction maps.
\end{interpretation}

A ubiquitous problem in sheaf theory is to decide whether $0$-cochains are coherent with respect to restriction maps. This consistency is expressed by the notion of a global section.

\begin{definition}[Global Sections]
    Suppose $\mathcal{F}$ is a cellular sheaf over $G$. A \define{global section} is a $0$-cochain $\mathbf{x} \in C^0(G; \mathcal{F})$ such that
    \begin{align}
        \mathcal{F}_{i \face ij}(x_i) = \mathcal{F}_{j \face ij}(x_j) \quad \forall i \in V, \forall j \in N_i.
        \label{eq:global-section}
    \end{align}
    The set of global sections is $\Gamma(G;\mathcal{F}) \subseteq C^0(G; \mathcal{F})$.
\end{definition}

\begin{interpretation}
    In a multi-agent coordination problem, agents are often required to satisfy locally-defined constraints that propagate globally, the simplest example being consensus. The global section condition (\cref{eq:global-section}) allows these local constraints to be encoded as restriction maps. For instance, the constant sheaf encodes a consensus constraint when $G$ is connected. This follows from $\Gamma(G;\underline{\mathbb{R}}^k) \cong \mathbb{R}^k$. If $G$ has $c$ connected components, it follows that $\Gamma(G;\underline{\mathbb{R}}^k) \cong \left(\mathbb{R}^k\right)^c$. A global section, then, corresponds to a particular consensus vector assigned to each connected component.
\end{interpretation}

\begin{table*}[htb]
\caption{Sheaf Theory Concepts and Potential Functions in Multi-Agent Systems}
\centering
\begin{tabular}{l c l c c l}
\toprule
\multicolumn{3}{c}{\textbf{Sheaf Theory Concepts}} & & \multicolumn{2}{c}{\textbf{Potential Functions}} \\
\midrule
\textbf{Concept} & \textbf{Notation} & \textbf{Interpretation} & & $U_e(y)$& \textbf{Coordination Goal} \\
\midrule
graph            & $G = (V, E)$      & bidirectional communication network & & $(1/2)\|y\|_2^2$ & Consensus \\
vertex stalk     & $\mathcal{F}(i)$  & state space of agent $i$ & & $y^\top Ay$ & Matrix-weighted consensus \cite{trinh_matrix-weighted_2018} \\
local section    & $x_i$             & local state & & $-(1/2)\|y\|_2^2$ & Dissensus \cite{franci_breaking_2023} \\
edge stalk       & $\mathcal{F}(ij)$ & communication space between agent $i$ and $j$ & & $(1/2)\|y-b\|_2^2$ & Reach displacement of $b$ \\
restriction map  & $\mathcal{F}_{i\face ij}$ & how agent $i$ sends messages to agent $j$ & & $(\|y\|_2^2-r^2)^2$ & Reach distance of $r$ \\
sheaf            & $\mathcal{F}$     & multi-agent communication structure & & & \\
0-cochains       & $C^0(G; \mathcal{F})$ & global state space & & & \\
0-cochain        & $\mathbf{x}$      & global state & & & \\
global sections  & $\Gamma(G; \mathcal{F})$ & feasible global states & & & \\
\bottomrule
\end{tabular}
\label{tab:combined}
\end{table*}

Deciding whether a $0$-cochain is a global section is not enough. Often, we want a flow to a global section from an arbitrary $0$-cochain. To construct such flows, we first introduce a few notions from homological algebra. Define the \define{coboundary operator}
$C^0(G; \mathcal{F}) \xrightarrow{\delta_{\mathcal{F}}} C^1(G; \mathcal{F})$ given by $(\delta_{\mathcal{F}} \mathbf{x})_{ij} = \mathcal{F}_{i \face ij}(x_i) - \mathcal{F}_{j \face ij}(x_j)$.
The \define{degree-0} and \define{degree-1 cohomology} of a cellular sheaf is given by 
\begin{align}
    \begin{aligned}
        H^0(G; \mathcal{F}) &= \ker \delta_{\mathcal{F}} \\
        H^1(G; \mathcal{F}) &= C^1(G; \mathcal{F})/ \image \delta_{\mathcal{F}}
    \end{aligned}
\end{align}
It follows immediately that $H^0(G; \mathcal{F})= \Gamma(G; \mathcal{F})$ because $x \in \ker \delta_{\mathcal{F}}$ precisely when $\mathcal{F}_{i \face ij}(x_i) - \mathcal{F}_{j \face ij}(x_j) = 0$ for all $ij \in E$. Global sections, then, can be computed as the kernel of the coboundary operator, where the relevant matrix grows in the number of nodes of the graph and the dimensions of the stalks.

A decentralized approach to computing global sections is to construct an appropriate local operator on $C^0(G;\mathcal{F})$ whose dynamics converge to a global section. A linear operator $L_\mathcal{F}$ -- a \define{linear sheaf Laplacian} -- whose construction relies on the Hodge Laplacian \cite{eckmann_harmonische_1944} leads to a heat equation $\dot{\mathbf{x}} = -L_{\mathcal{F}} \mathbf{x}$ converging to a global section for an arbitrary initial condition \cite[Proposition 8.1]{hansen_toward_2019}. To analyze a wider class of coordination problems, we utilize a generalization of the linear sheaf Laplacian \cite[\S 10]{hansen_opinion_2021}. In addition to the specification of a graph and a sheaf, the nonlinear sheaf Laplacian is determined by a choice of edge potentials.

\begin{definition}[Nonlinear Sheaf Laplacian]
    Suppose $\mathcal{F}$ is a cellular sheaf over $G$, and suppose $\{U_e: \mathcal{F}(e) \to \mathbb{R} \}_{e \in E}$ are (possibly nonlinear) potential functions. Then, the \define{nonlinear sheaf Laplacian} is the map $ C^0(G;\mathcal{F}) \xrightarrow{L^{\nabla U}_\mathcal{F}} C^0(G;\mathcal{F})$
    defined $L^{\nabla U}_\mathcal{F} = \delta^\top \circ \nabla U \circ \delta$
    where $\nabla U: C^1(G;\mathcal{F}) \to C^1(G;\mathcal{F})$ is defined as the gradient of $U(\mathbf{y}) = \sum_{e \in E} U_e(y_e)$. 
\end{definition}

\begin{examples}
     Given $G$ and $\mathcal{F}$, suppose $U_{e}(y_e) = (1/2) \|y_e\|^2_2$ for all $e \in E$. Then, $L_\mathcal{F}^{\nabla U} = \delta^\top \circ \delta$ and is denoted $L_\mathcal{F}$. In particular, the linear sheaf Laplacian $L_{\underline{\mathbb{R}}^k}$ of the constant sheaf on $G$ is the Kronecker product $L_G \otimes I_{k \times k}$ where $L_G$ is the unnormalized graph Laplacian $L_G$, or simply $L_G$ in the one-dimensional case.
\end{examples}

The nonlinear sheaf Laplacian acts on local sections as
    \begin{align*}
        \begin{aligned}
            (L_\mathcal{F}^{\nabla U} \mathbf{x})_i = \sum_{j \in N_i} \mathcal{F}_{i \face ij}^\top \biggl( \nabla U_{ij}  \bigl( \mathcal{F}_{i \face ij}(x_i) - \mathcal{F}_{j \face ij}(x_j) \bigr) \biggr).
        \end{aligned}
    \end{align*}
Consequently, the nonlinear sheaf Laplacian can be computed locally in the network because $(L_\mathcal{F}^{\nabla U} \mathbf{x})_i$ only requires information about the local sections, restriction maps and of neighbors $N_i$, and potential functions $\{U_{ij}\}_{j \in N_i}$.

While linear sheaf Laplacians converge to global sections specifying local linear constraints, nonlinear Laplacian dynamics are more subtle. For example, nonconvex potentials can lead to diverging trajectories. For instance,
given $G$ and $\mathcal{F}$, suppose $U_{e}(y) = -(1/2) \|y\|^2_2$ for all $e \in E$. Then, $L_\mathcal{F}^{\nabla U} = -L_\mathcal{F}$, and trajectories with $\mathbf{x}(0) \notin H^0(G; \mathcal{F})$ diverge from their the projection onto $H^0(G;\mathcal{F})$.

%%%%%%%%%%%%%%%%%%%%%%%%%%%%%%%%%%%%%%%%%%%%%%%%%%%%%%%%%%%%%%%%%%%%%%%%%%
\section{Nonlinear Homological Programming for Multi-agent Coordination}
\label{sec:homological-programming}
%%%%%%%%%%%%%%%%%%%%%%%%%%%%%%%%%%%%%%%%%%%%%%%%%%%%%%%%%%%%%%%%%%%%%%%%%%

This section presents the primary modeling tool used in this framework -- nonlinear homological programming -- and formulates a general multi-agent coordination problem within this framework.

\subsection{Nonlinear Homological Programming}

The following definition extends homological programming~\cite{hansen_distributed_2019} with nonlinear constraints and objective functions valued in the extended real numbers $\Rbar = \mathbb{R} \cup \{\infty\}$.

\begin{definition}[Nonlinear Homological Program]
    A \define{nonlinear homological program} consists of the following data:
    \begin{enumerate}[leftmargin=*]
        \item Undirected graph $G=(V,E)$;
        \item Cellular sheaf $\mathcal{F}$ on $G$;
        \item Objective functions $ \{ f_i\maps \mathcal{F}(v)\to\Rbar \}_{i \in V}$;
        \item Potential functions $\{ U_{e}\maps\mathcal{F}(e)\to\mathbb{R} \}_{e \in E}$.
    \end{enumerate}
    The homological program defined by this data is the optimization problem
\begin{equation}
    \begin{array}{ll@{}ll}
    \minimize\limits_{\mathbf{x} \in C^0(G;\mathcal{F})} & \sum_{i\in V}f_i(x_i) \\
    \text{subject to} & L_\mathcal{F}^{\nabla U} \mathbf{x}= 0
    \end{array}
    \tag{P}
    \label{eq:nonlinear-homological-programming}
\end{equation}
\end{definition}

\begin{example}[Distributed Optimization]
    In a system of $N$ agents where each agent has an individual convex objective function $f_i(x)$ with decision variable $x \in \mathbb{R}^k$, a ubiquitous distributed optimization problem \cite{yang_survey_2019} is to minimize the sum of each agent's cost function: $\minimize_{x \in \mathbb{R}^k} \sum_{i = 1}^N f_i(x)$. A typical distributed solution strategy is to introduce local copies of $x$ for each agent to decouple the objective and use a coupling consensus constraint. This can be defined as a homological program for the constant sheaf $\underline{\R}^k$ on a connected agent communication topology, as shown in \cite{hansen_distributed_2019}.
\end{example}

The following theorem describes when nonlinear homological programs are convex optimization problems, enabling us to solve them efficiently.

\begin{theorem}\label{thm:hp-convexity}
    Let $\mathsf{P} = \left(V, E, \mathcal{F}, \{f_i\}, \{U_{e}\} \right)$ be a nonlinear homological program. If $U_e$ is differentiable and convex for each $e\in E$, and $f_i$ is convex for every $i \in V$, then $\mathsf{P}$ is a convex optimization problem.
    \footnote{For proofs of all theorems, please consult our technical report: \\ \url{https://arxiv.org/abs/2504.02049}}
\end{theorem}

For illustrative purposes, we introduce single-agent optimal control as a simple example of a homological program which will be further utilized when posing multi-agent optimal control problems.

\begin{example}[Single-Agent Optimal Control]\label{ex:single-agent}
Single-agent optimal control problems for a fixed time horizon $T$ form homological programs on cellular sheaves whose underlying graphs are $T+1$ length paths. To see why, consider a control system with states in $\R^n$, control inputs in $\R^m$, and discrete dynamics $x(t+1) = Ax(t) + Bu(t)$.
We model this agent's dynamics for $T$ time-steps from an initial condition $x(1)=c$ as a cellular sheaf $\mathcal{D}$ on the $(T+1)$-vertex path graph $P_{T+1}$ (see \cref{fig:LTI-sheaf}). The overall potential function for $\mathcal{D}$ is $\Phi(\mathbf{y}) = (1/2)\|y_0-c\|_2^2 + \sum_{e \in E} (1/2)\|y_e\|_2^2$ for all $\mathbf{y}\in C^1(P_{T+1};\mathcal{D})$. Zeros of the associated nonlinear Laplacian $L_\mathcal{D}^{\nabla \Phi}=\delta_\mathcal{D}^\top\circ\nabla\Phi\circ\delta_\mathcal{D}$ then correspond precisely to admissible trajectories of the system across $T$ time-steps starting from $x(1)=c$.

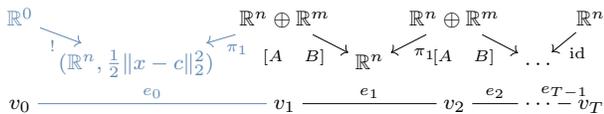
\begin{figure}[htb]
 \footnotesize
% https://q.uiver.app/#q=WzAsMTIsWzIsMCwiXFxSXntufSBcXG9wbHVzXFxSXnttfSJdLFszLDEsIlxcUl5uIl0sWzQsMCwiXFxSXm5cXG9wbHVzXFxSXm0iXSxbNSwxLCJcXGNkb3RzIl0sWzYsMCwiXFxSXm4iXSxbMiwyLCJ2XzEiXSxbNCwyLCJ2XzIiXSxbNiwyLCJ2X1QiXSxbNSwyLCJcXGNkb3RzIl0sWzEsMSwiKFxcUl5uLCBcXGZyYWN7MX17Mn1cXHx4LWNcXHxfMl4yKSJdLFswLDAsIlxcUl4wIl0sWzAsMiwidl8wIl0sWzAsMSwiW0EgXFxxdWFkIEJdIiwyXSxbMiwxLCJcXHBpXzEiXSxbMiwzLCJbQVxccXVhZCBCXSIsMl0sWzQsMywiXFxtYXRocm17aWR9Il0sWzUsNiwiZV8xIiwwLHsic3R5bGUiOnsiaGVhZCI6eyJuYW1lIjoibm9uZSJ9fX1dLFs2LDgsImVfMiIsMCx7InN0eWxlIjp7ImhlYWQiOnsibmFtZSI6Im5vbmUifX19XSxbNyw4LCJlX3tULTF9IiwyLHsic3R5bGUiOnsiaGVhZCI6eyJuYW1lIjoibm9uZSJ9fX1dLFswLDksIlxccGlfMSIsMCx7ImNvbG91ciI6WzAsNjAsNjBdfSxbMCw2MCw2MCwxXV0sWzEwLDksIiEiLDIseyJjb2xvdXIiOlswLDYwLDYwXX0sWzAsNjAsNjAsMV1dLFsxMSw1LCJlXzAiLDAseyJjb2xvdXIiOlswLDYwLDYwXSwic3R5bGUiOnsiaGVhZCI6eyJuYW1lIjoibm9uZSJ9fX0sWzAsNjAsNjAsMV1dXQ==
\[\begin{tikzcd}[cramped,sep=tiny] 
	{\color{gatas-blue}\R^0} && {\R^{n} \oplus\R^{m}} && {\R^n\oplus\R^m} && {\R^n} \\
	& {\color{gatas-blue}(\R^n, \frac{1}{2}\|x-c\|_2^2)} && {\R^n} && \cdots \\
	{v_0} && {v_1} && {v_2} & \cdots & {v_T}
	\arrow["{!}"', color={rgb,255:red,97;green,136;blue,178}, from=1-1, to=2-2]
	\arrow["{\pi_1}", color={rgb,255:red,97;green,136;blue,178}, from=1-3, to=2-2]
	\arrow["{[A \quad B]}"', from=1-3, to=2-4]
	\arrow["{\pi_1}", from=1-5, to=2-4]
	\arrow["{[A\quad B]}"', from=1-5, to=2-6]
	\arrow["{\mathrm{id}}", from=1-7, to=2-6]
	\arrow["{e_0}", color={rgb,255:red,97;green,136;blue,178}, no head, from=3-1, to=3-3]
	\arrow["{e_1}", no head, from=3-3, to=3-5]
	\arrow["{e_2}", no head, from=3-5, to=3-6]
	\arrow["{e_{T-1}}"', no head, from=3-7, to=3-6]
\end{tikzcd}\]
\caption{The cellular sheaf encoding the evolution of a discrete LTI system from an initial condition $x(1)=c$ over $T$ time-steps. Note that $!\maps \R^0\to\R^n$ is the unique map from $\R^0$, and $\pi_1$ is the first projection. The graph which the sheaf is defined on is given on the bottom with the stalks and restriction maps lying above. The blue section uses an edge potential function to fix the initial condition to $x(1)=c$. All other edge potential functions are the standard consensus potentials so we omit them from the diagram for clarity. The black section then encodes the evolution of the dynamics from $c$ for $T$ time-steps.}
\label{fig:LTI-sheaf}
\end{figure}

Now, suppose we wish to drive the system to a desired state using stage cost functions $f^t\maps \R^n\times \R^m\to\R$ for $t\in[T-1]$. A typical multistage optimal control problem starting from initial state $x(1)=c$ is
\begin{equation}\label{eq:opt-control}
\begin{array}{ll@{}ll}
    \text{minimize} & \sum_{t=1}^{T-1}f^t(x(t),u(t)) \\
    \text{subject to} & x(t+1) = Ax(t)+Bu(t)\quad\forall t\in [T-1]\\
    & x(1) = c.
\end{array}
\end{equation}
Note that our framework can easily incorporate other convex constraints on the states and controls; however, we leave these out for simplicity of presentation.
For notational convenience, let $\mathbf{x}=\text{vec}(x(1),\dots,x(T))$ and $\mathbf{u}=\text{vec}(u(1),\dots,u(T-1))$, where vec denotes the concatenation of column vectors into a new column vector. We have already seen that the constraints in \cref{eq:opt-control} correspond to the zeros of the nonlinear Laplacian for $\mathcal{D}$. Thus the problem in \cref{eq:opt-control} is equivalent to the homological program
\begin{equation}\label{eq:opt-control-hp}
\begin{array}{ll@{}ll}
    \text{minimize} & \sum_{t=1}^{T-1}f^t(x(t),u(t)) \\
    \text{subject to} & L_{\mathcal{D}}^{\nabla \Phi}(\mathbf{x},\mathbf{u})=0.
\end{array}
\end{equation}    
\end{example}

%---------------------------------------------------
\subsection{General Multi-Agent Coordination Homological Program}
%---------------------------------------------------

Now suppose we have a multi-agent control system with $N$ agents on a communication topology $G=(V,E)$. Each agent $i\in[N]$ has state space $\R^{n_i}$, control space $\R^{m_i}$, and discrete dynamics $x_i(t+1)=A_ix_i(t) + B_iu_i(t)$. The total system dynamics are then given by the direct sum of the dynamics for each agent. Let $Q = \sum_{i\in[N]} n_i$ and $R=\sum_{i\in[N]}m_i$ denote the global state and control input dimensions respectively. For notational convenience, we introduce block matrix variables $\mathbf{X}\in \R^{Q\times T}$ and $\mathbf{U}\in\R^{R\times T}$ to represent $T$-length trajectories of the global state space (resp. control space), where block $i,t$ is in $\R^{n_i\times 1}$ (resp. $\R^{m_i\times 1}$). Then, $\mathbf{X}[:,t]$ denotes the global state at time $t$, $\mathbf{X}[i,:]$ denotes the state trajectory of agent $i$ across all time-steps, and $\mathbf{X}[i,t]$ denotes the state of agent $i$ at time-step $t$. The same goes for $\mathbf{U}$.

A general distributed multi-agent coordination problem is that of choosing control inputs to drive this global system to a desired goal while utilizing only local computations and communication between agents. We formulate this as a multi-stage optimal control problem over a fixed time horizon $T$. We allow each agent to have its own local optimal control problem of the same form as \cref{eq:opt-control}. Specifically, each agent $i\in[N]$ has stage costs $f_i^t$ and dynamic evolution sheaf $\mathcal{D}_i$ with potential function $\Phi_i$ from a given initial condition $x_i(1)=c_i$, as in \cref{ex:single-agent}. Let $\mathcal{C}_i$ denote the set $\{ (\mathbf{x},\mathbf{u})\in C^0(P_{T+1};\mathcal{D}_i)\mid L_{\mathcal{D}_i}^{\nabla \Phi_i}(\mathbf{x},\mathbf{u})=0\}$. Then we define the functional form of \cref{eq:opt-control-hp} as $J_i\maps C^0(P_{T+1},\mathcal{D}_i)\to \Rbar$ given by
\begin{equation}\label{eq:cost-to-go}
    J_i(\mathbf{x},\mathbf{u}) = \sum_{t=1}^{T-1}f_i^t(x(t),u(t)) + \chi_{\mathcal{C}_i}(\mathbf{x},\mathbf{u}),
\end{equation}
where $\chi_{\mathcal{C}_i}$ denotes the convex indicator function of the set $\mathcal{C}_i$. The global objective function $\sum_{i=1}^NJ_i(\mathbf{X}[i,:],\mathbf{U}[i,:])$ is then a decoupled multi-agent optimal control problem where each agent wants to fulfill its individual control goal disregarding other agents.

To coordinate our multi-agent system, we introduce a cellular sheaf $\mathcal{F}$ on $G$ with the interpretations given in Table~\ref{tab:combined} (left). We can then formulate the coordination goal between each pair of communicating agents $e=i\sim j\in E$ by choosing an appropriate edge potential function $U_e$ from Table~\ref{tab:combined} (right). Such a choice of potential functions determines a global potential function $U(\mathbf{y})=\sum_{e\in E}U_e(\mathbf{y}_e)$ defined over 1-cochains and yields the coordination objective
\[
\underset{\mathbf{X}[\colon,t]\in C^0(G;\FF)}{\textnormal{minimize }} U(\delta_\FF \mathbf{X}[\colon,t])
\]
for a given $t\in [T]$.
Taken together, this lets us define a general multi-agent coordination problem as
\begin{equation}\label{eq:prob-relaxation}
    \textnormal{minimize }\sum_{i=1}^N J_i(\mathbf{X}[i,:],\mathbf{U}[i,:]) + \gamma U(\delta_\FF\mathbf{X}[:,T]),
\end{equation}
In other words, this objective asks each agent in the MAS to achieve its own control objective as much as possible subject to its dynamic constraints and initial condition while also cooperating with the other agents to achieve the coordination goal by the end of the time horizon. The free parameter $\gamma$ can be increased to emphasize the importance of the group objective over the individual objectives. When $U$ is convex and bounded below, minimizing $U\circ\delta_\FF$ is equivalent to finding a point where $\nabla (U\circ\delta_\FF) = \delta_\FF^\top\circ \nabla U\circ\delta_\FF = 0$. This gives us a strengthened multi-agent homological program:
\begin{equation}\label{eq:general-prob}
\begin{array}{ll@{}ll}
     \underset{\mathbf{X},\mathbf{U}}{\text{minimize}}& \sum_{i=1}^N J_i(\mathbf{X}[i,:],\mathbf{U}[i,:])  \\
     \text{subject to}& L_\mathcal{F}^{\nabla U}\mathbf{X}[:,T]=0,
\end{array}
\end{equation}
where the coordination goal is imposed as a hard constraint rather than an objective term.

\begin{remark}
    It should be noted that there are multiple sheaf homological constraints present in \cref{eq:general-prob}. Specifically, the sheaves $\mathcal{D}_i$ encode the dynamic evolution of each agent $i\in[N]$ from an initial condition to the end of the time horizon. Meanwhile, the sheaf $\mathcal{F}$ and its associated nonlinear Laplacian specifies the coordination constraint which must be satisfied by the end of the time horizon. We can think of this as a sort of \emph{nested} or \emph{hierarchical} homological program, where the objectives on each node of the coordination sheaf are themselves homological programs.
\end{remark}

%%%%%%%%%%%%%%%%%%%%%%%%%%%%%%%%%%%%%%%%%%%%%%%%%%
\section{Solving Homological Programs with ADMM}
\label{sec:ADMM}
%%%%%%%%%%%%%%%%%%%%%%%%%%%%%%%%%%%%%%%%%%%%%%%%%%

We examine how we can use ADMM to derive a distributed solver for homological programs, thus allowing the design of distributed controllers for achieving multi-agent coordination.
Algorithm 1 presents a distributed solution algorithm for a general nonlinear homological program $\mathsf{P}= (V,E,\mathcal{F},\{f_i\},\{U_e\})$ subject to the following assumptions.

\begin{assumption}\label{as:convex-p}
The homological program $\mathsf{P}$ satisfies the convexity conditions of Theorem \ref{thm:hp-convexity}, and each node objective $f_i$ is closed and proper.
\end{assumption}

\begin{assumption}\label{as:quad-potentials}
    Each edge potential $U_e$ is strongly convex with unique minimizer of $b_e$ such that $U_e(b_e)=0$.
\end{assumption}

To use ADMM to solve $\mathsf{P}$, we first introduce extra variables $z_i\in \mathcal{F}(i)$ for each vertex stalk and rewrite the problem in the equivalent consensus form:
\begin{equation}\label{eq:consensus-form}
\begin{array}{ll@{}ll}
    \textnormal{minimize } & \sum_i f_i(x_i) + \chi_\mathcal{C}(\mathbf{z})  \\
    \textnormal{subject to } & \mathbf{x}-\mathbf{z} = 0 
\end{array}
\end{equation}
where $\mathcal{C}=\{\mathbf{z}\in C^0(G;\mathcal{F})\mid L^{\nabla U}_\mathcal{F}\mathbf{z}=0\}$. This problem has the following augmented Lagrangian with scaled dual variable as in \cite{boyd_distributed_2010}:
\[
\mathcal{L}_\rho(\mathbf{x},\mathbf{z},\mathbf{y}) = \sum_i f_i(x_i) + \chi_\mathcal{C}(\mathbf{z}) + (\rho/2)\|\mathbf{x} -\mathbf{z} + \mathbf{y}\|_2^2.
\]
Applying ADMM gives the iterative update rule:
\begin{equation}\label{eq:admm}
\begin{array}{ll@{}ll}
    x_i^{k+1} & \coloneqq \argmin_{x_i} f_i(x_i) + (\rho/2)\|x_i - z_i^k + y_i^k\|_2^2  \\
    \mathbf{z}^{k+1} & \coloneqq \Pi_\mathcal{C}(\mathbf{x}^{k+1} + \mathbf{y}^k) \\
    y_i^{k+1} & \coloneqq y_i^k + x_i^{k+1} - z_i^{k+1} 
\end{array}
\end{equation}
where $\Pi_\mathcal{C}$ denotes Euclidean projection onto $\mathcal{C}$.
The dynamics of the $x_i$ and $y_i$ updates can be computed in parallel. The core difficulty of the algorithm then is computing the $\mathbf{z}$ update in a distributed fashion. It turns out that running the nonlinear heat equation, which is a local operation, converges precisely to the desired projection. We prove this in the following theorem.

\begin{theorem}\label{thm:nonlinear-projection}
    Let $\mathcal{F}$ be a cellular sheaf on a graph $G=(V,E)$. For each edge $e\in E$, let $U_e\maps \mathcal{F}(e)\to \R$ be a differentiable, strongly convex function with unique global minimum point at $b_e$ with minimum value $U_e(b_e)=0$. Define $U=\sum_{e\in E}U_e$ and $\mathbf{b} = \textnormal{vec}(b_1,\dots,b_{|E|})$. If $\mathbf{b}\in \image \delta_\mathcal{F}$, then trajectories of the nonlinear Laplacian dynamics
    \begin{equation}\label{eq:nonlinear-diffusion-dynamics}
    \dot{\mathbf{x}} = -\alpha L_\mathcal{F}^{\nabla U}\mathbf{x}
    \end{equation}
    converge to the orthogonal projection of the initial condition onto $\delta_\FF^+\mathbf{b} + H^0(G;\FF) =\ker L_\FF^{\nabla U}$ for a given diffusivity $\alpha>0$, where $\delta_\FF^+$ denotes the Moore-Penrose pseudoinverse.
    
    %the affine subspace $\{x\mid \delta x -b=0\}$ of $C^0(G;\mathcal{F})$.
\end{theorem}

This gives us Algorithm \ref{alg:admm} for solving nonlinear homological programs.
\SetKwComment{Comment}{/* }{ */}
\begin{algorithm}[htbp]
\footnotesize
    \caption{Distributed Solve}\label{alg:admm}
    \KwIn{Sheaf $\FF$, potential functions $U_e$, and initial state 
    $(\mathbf{z}^0,\mathbf{y}^0)$}
    \KwIn{Step size $\rho$, diffusivity $\alpha$, tolerances $\epsilon_1,\epsilon_2$, and $K\in \mathbb{N}$}
    \KwOut{$\mathbf{z}^*$,$\mathbf{y}^*$}
    $\mathbf{z}\gets\mathbf{z}^0$\;
    $\mathbf{y}\gets\mathbf{y}^0$\;
    \While{$j \leq K$}{
        \For{$i\gets 1$ \KwTo $N$}{
            $x_i \gets \argmin_{x_i} f_i(x_i) + (\rho/2)\|x_i - z_i + y_i\|^2_2$\;
        }
        $\mathbf{x}\gets \textnormal{vec}(\{x_i\}_{i\in[N]})$\;
        $\mathbf{z}\gets \text{sheafDiffusion}(\FF, U, \alpha, \epsilon_2, \mathbf{x} +\mathbf{y})$\;
        \If{$\mathbf{x}-\mathbf{z} < \epsilon_1$}{
            break\;
        }
        \For{$i\gets 1$ \KwTo $N$}{
            $y_i\gets y_i + x_i - z_i$\;
        }
        
    }
    \Return{$\mathbf{z},\mathbf{y}$}
\end{algorithm}
In this algorithm, each iteration of the two for-loops can be computed in parallel by each agent. The operation $\textnormal{sheafDiffusion}(\FF, U, \alpha, \epsilon, \mathbf{z}^0)$ denotes running the nonlinear heat equations for $\FF$ with diffusivity $\alpha$ to $\epsilon$ convergence from initial condition $\mathbf{z}^0$.

\begin{theorem} \label{thm:admm}
    Suppose a homological program $\mathsf{P}$ satisfies assumptions \ref{as:convex-p} and \ref{as:quad-potentials}. Assume also that the Lagrangian of \cref{eq:consensus-form} has a saddle point. Then Algorithm \ref{alg:admm} applied to $\mathsf{P}$ has the following properties:
    \begin{itemize}
        \item Residual convergence: $\mathbf{x}^k-\mathbf{z}^k\to 0$ as $k\to\infty$
        \item Objective convergence: $\sum_i f_i(x_i) + \chi_\mathcal{C}(\mathbf{z})\to p^*$ as $k\to \infty$.
        \item Dual variable convergence: $\mathbf{y}^k\to\mathbf{y}^*$ as $k\to \infty$ where $\mathbf{y}^*$ is a dual optimal point.
    \end{itemize}
\end{theorem}

\begin{remark}
    When potential functions are chosen that violate Assumption \ref{as:quad-potentials}, the same ADMM algorithm can still be used to solve the relaxed problem \cref{eq:prob-relaxation} with slight modification to the $\mathbf{z}$ update.
\end{remark}

%%%%%%%%%%%%%%%%%%%%%%%%%%%%%%%%%%%%%%%%%%%%%%%%%%%%%%%%%%%%%%
\section{Case Study: Multi-Domain Team Vehicle Coordination}
\label{sec:examples}
%%%%%%%%%%%%%%%%%%%%%%%%%%%%%%%%%%%%%%%%%%%%%%%%%%%%%%%%%%%%%%

Multi-vehicle systems (MVS) provide a rich class of examples of multi-agent coordination problems that can be formulated and solved as homological programs. In this section, we describe some of these coordination problems in more detail through a case study and explain how they are instances of homological programs.

\begin{figure}
    \centering
    \includegraphics[width=0.75\linewidth]{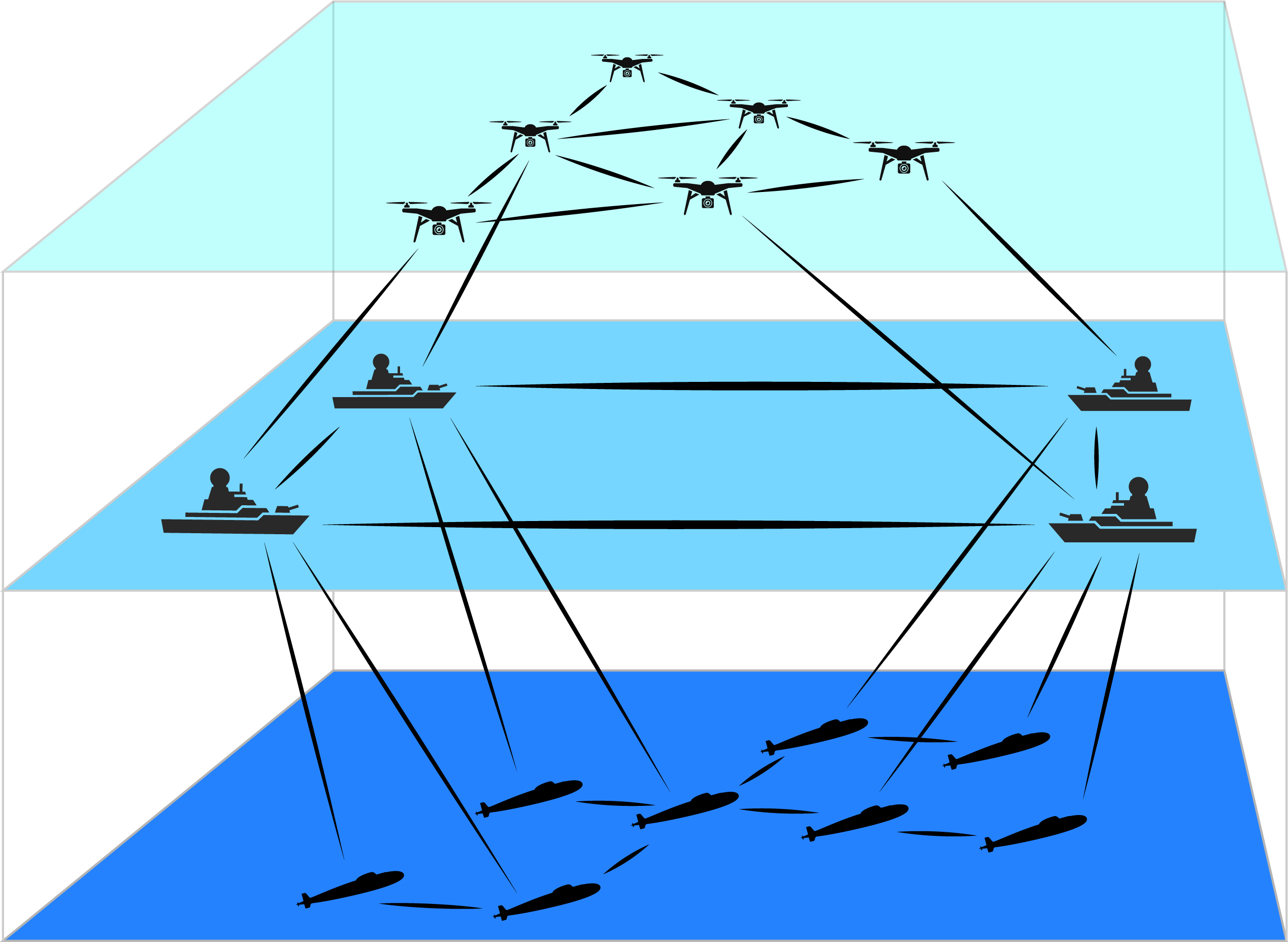}
    \caption{Multi-domain operation: UAVs, USVs, \& UUVs.}
    \label{fig:multi-domain}
\end{figure}

We consider a team of heterogeneous autonomous vehicles (see \cref{fig:multi-domain}) that includes unmanned aerial vehicles (UAV), unmanned surface vehicles (USV), and unmanned underwater vehicles (UUV). The vehicles are indexed by $V = \{1,2,\dots, N\}$. In this setting, communication modalities (e.g.~radio frequency, sonar, optical, or tethered links) are domain-dependent, varying with the operational environment of each vehicle type.

We assume vehicles have a fixed communication network. UAVs have a dense communication network between other UAVs, and UAVs also communicate with USVs. USVs serve as relays for communication between USVs and UAVs. Consequently, USVs have links with other UAVs as well other USVs and links with designated UUVs. Finally, UUVs communicate with other UUVs, and UUVs communicate with designated USVs. We model the entire communication network according to an undirected graph $G = (V,E)$ where an edge $e \in E$ denotes the presence of a bidirectional communication channel between the vehicles at its endpoints. Let $G_{\mathrm{UAV}}$, $G_{\mathrm{USV}}$, and 
$G_{\mathrm{UUV}}$ denote the respective subgraphs corresponding to each respective inter-team communication network, and let $G_{\mathrm{UAV-USV}}$ and $G_{\mathrm{UUV-USV}}$ denote the subgraphs corresponding to each intra-team network.

Vehicles navigate in a rectangular domain $D \subseteq \mathbb{R}^3$, which includes airspace, a surface region, and an undersea region. We assume 
that for each~$i \in V$ the dynamics of 
vehicle~$i$ are modeled by double integrator dynamics
\begin{align}
\begin{aligned}
    \dot{p}_i &= v_i \\
    \dot{v}_i &= u_i, \label{eq:double-integrator}
\end{aligned}
\end{align}
where $p_i \in \mathbb{R}^d$ is agent~$i$'s position, $v_i \in \mathbb{R}^d$ is agent~$i$'s velocity, and $u_i \in \mathbb{R}^d$ is agent~$i$'s control input. Then, the state of agent~$i$ is~$x_i = \text{vec}(p_i,v_i)
%\begin{bmatrix} p_i \\ v_i\end{bmatrix}
$. 
The dimensions of the agents' state spaces depend on the type of vehicle: UUVs and UAVs move in~$\R^3$ and hence have~$d=3$, 
while~USVs are modeled as moving in~$\R^2$ and hence have~$d=2$.
We assume that agents fully observe their own position and velocity.
Consider the following scenario which we describe as an optimal control problem: 
\begin{itemize}[leftmargin=*]
    \item The USV team is tasked with making a translationally-invariant grid formation that is specified by displacement vectors $\hat{p}_{i,j}$. A subset of edges $E_{\mathrm{USV}} \subseteq E_{\mathrm{surface}}$ indicates which pairs of agents must maintain a specified relative position in order to assemble
    the given formation. 
    \item The UAV team flocks together, forming consensus on their speeds and bearings while maintaining a fixed altitude $h_{\mathrm{UAV}}$. The UAVs also need to remain within a given distance $r_{\mathrm{UAV-USV}}$ from the USVs that they communicate with in order to remain within communication range. 
    \item The UUV team flocks together, forming consensus on their speeds and bearings, but are free to vary their depth from the surface. In order to maintain communication, UUVs also need to remain within a given distance $r_{\mathrm{UUV-USV}}$ to USVs,
    and they also need to remain with a given distance of 
    other UUVs for the same reason. 
    \item All vehicles minimize the actuation required for their tasks.
\end{itemize}
With the interpretation (see \cref{sec:prelim}) of stalks as local state-spaces, $0$-cochains as global state-spaces, global sections as global coordination constraints, and edge potentials as agent-to-agent coordination objectives, we synthesize a homological program $\mathsf{P} = (V,E,\mathcal{F}, \{U_e\}, \{f_i\})$ describing the above complex coordination problem.

\paragraph{Stalks}
The cellular sheaf $\mathcal{F}$ is defined over the communication graph $G=(V,E)$. For each node $i \in V$, the vertex stalk is assigned based on the vehicle type:
{\footnotesize \begin{align*}
    \mathcal{F}(i) =
    \begin{cases}
    \mathbb{R}^3 \oplus \mathbb{R}^3 & \mathrm{if}~i \in V_{\mathrm{UAV}} \cup V_{\mathrm{UUV}}  \\
    \mathbb{R}^2 \oplus \mathbb{R}^2 & \mathrm{if}~i \in V_{\mathrm{USV}}
    \end{cases}.
\end{align*}}
For each edge $e = ij \in E$, the edge stalk is defined based on the type of connection. 
The edge stalks are:
{\footnotesize \begin{align*}
    \mathcal{F}(e) = 
    \begin{cases}
        \mathbb{R}^3 \oplus \mathbb{R}^3 & \mathrm{if}~e \in E_{\mathrm{UAV}} \cup E_{\mathrm{UUV}} \\
        \mathbb{R}^2 & \mathrm{if}~e \in E_{\mathrm{USV}} \\
        \mathbb{R}^3 & \mathrm{if}~e \in E_{\mathrm{UAV-USV}} \\
        \mathbb{R}^3 & \mathrm{if}~e \in E_{\mathrm{UUV-USV}}
    \end{cases}
\end{align*}}

\subsubsection*{Restriction maps}
In this case study, restriction maps serve as an interface between agent types as well as communication modalities. Thus, restriction maps $\mathcal{F}_{i \face e}: \mathcal{F}(i) \to \mathcal{F}(e)$ project the state of each agent to the shared edge state-space:
{\footnotesize \begin{itemize}[leftmargin=*]
    \item For $e = ij \in E_{\text{UAV}}$:
        \[
        \mathcal{F}_{i \face e}(x_i) = \begin{bmatrix} p_i \\ v_i \end{bmatrix}, \quad \mathcal{F}_{j \face e}(x_j) = \begin{bmatrix} p_j \\ v_j \end{bmatrix}
        \]
        so $y_e =  (\delta \mathbf{x})_e = \begin{bmatrix} p_i - p_j \\ v_i - v_j \end{bmatrix} \in \mathbb{R}^3 \oplus \mathbb{R}^3$.
    \item For $e = ij \in E_{\text{UUV}}$:
        \[
        \mathcal{F}_{i \face e}(x_i) = \begin{bmatrix} p_i \\ v_i \end{bmatrix}, \quad \mathcal{F}_{j \face e}(x_j) = \begin{bmatrix} p_j \\ v_j \end{bmatrix}
        \]
        so $y_e = (\delta \mathbf{x})_3 = \begin{bmatrix} p_i - p_j \\ v_i - v_j \end{bmatrix} \in \mathbb{R}^3 \oplus \mathbb{R}^3$.
    \item For $e = ij \in E_{\text{USV}}$:
        \[
        \mathcal{F}_{i \face e}(x_i) = \begin{bmatrix} p_{i,x} \\ p_{i,y} \end{bmatrix}, \quad \mathcal{F}_{j \face e}(x_j) = \begin{bmatrix} p_{j,x} \\ p_{j,y} \end{bmatrix}
        \]
        so $y_e = (\delta \mathbf{x})_e = \begin{bmatrix} p_{i,x} - p_{j,x} \\ p_{i,y} - p_{j,y} \end{bmatrix} \in \mathbb{R}^2$.
    \item For $e = ij \in E_{\text{UAV-USV}}$:
        \[
        \mathcal{F}_{i \face e}(x_i) = p_i, \quad \mathcal{F}_{j \face e}(x_j) = \begin{bmatrix} p_{j,x} \\ p_{j,y} \\ 0 \end{bmatrix}
        \]
        so $y_e = (\delta \mathbf{x})_e = p_i - \begin{bmatrix} p_{j,x} \\ p_{j,y} \\ 0 \end{bmatrix} \in \mathbb{R}^3$.
    \item For $e = ij \in E_{\text{UUV-USV}}$:
        \[
        \mathcal{F}_{i \face e}(x_i) = p_i, \quad \mathcal{F}_{j \face e}(x_j) = \begin{bmatrix} p_{j,x} \\ p_{j,y} \\ 0 \end{bmatrix}
        \]
        so $y_e = (\delta \mathbf{x})_e = p_i - \begin{bmatrix} p_{j,x} \\ p_{j,y} \\ 0 \end{bmatrix} \in \mathbb{R}^3$.
\end{itemize}}

\subsubsection*{Edge Potentials}
Edge potentials encode heterogeneous coordination objectives. Edge potentials $\{U_e\}$ are defined as follows:
{\footnotesize \begin{itemize}[leftmargin=*]
    \item For $e = ij \in E_{\text{form}} \subseteq E_{\text{USV}}$:
        \[
        U_e(y_e) = \frac{1}{2} \| y_e - \hat{p}_{i,j} \|_2^2
        \]
        where $\hat{p}_{i,j} \in \mathbb{R}^2$ is the desired displacement. For edges $e \in E_{\text{USV}} \setminus E_{\text{form}}$ not involved in formation, we set $U_e(y_e) = 0$.
    \item For $e = ij \in E_{\text{UAV}}$:
        \[
            U_e(y_e) = \frac{1}{2} \| v_i - v_j \|_2^2.
        \]
    \item For $e = ij \in E_{\text{UUV}}$:
        \[
        U_e(y_e) = \frac{1}{2} \| v_i - v_j \|_2^2 + \frac{1}{2} ( \| p_i - p_j \|_2^2 - r_{\text{UUV}}^2 )^2
        \]
        where $r_{\text{UUV}}>0$ is the desired UUV-UUV communication distance.
    \item For $e = ij \in E_{\text{UAV-USV}}$:
        \[
        U_e(y_e) = \frac{1}{2} ( \| y_e \|_2 - r_{\text{UAV-USV}} )^2
        \]
        where $r_{\text{UAV-USV}}$ is the desired UAV-USV communication distance.
    \item For $e = ij \in E_{\text{UUV-USV}}$:
        \[
        U_e(y_e) = \frac{1}{2} ( \| y_e \|_2^2 - r_{\text{UUV-USV}^2} )^2
        \]
        where $r_{\text{UUV-USV}}$ is the desired UUV-USV communication distance.
\end{itemize}}

\subsubsection*{Objective Functions}
For each agent $i \in V$, the objective is to minimize the control effort over a time horizon $T$, subject to dynamics and coordination constraints:
{\footnotesize \begin{align*}
f_i(x_i, u_i) = \sum_{t=1}^{T} \| u_i^t \|_2^2
\end{align*}}
where $u_i^t \in \mathbb{R}^3$ (for UAVs and UUVs) or $\mathbb{R}^2$ (for USVs) is the control input at time $t$, and the optimization is subject to the dynamics (\cref{eq:double-integrator}), and the coordination constraints encoded by the sheaf $\mathcal{F}$ and potentials $\{U_e\}$.

%%%%%%%%%%%%%%%%%%%%%%%%%%%%%%%%%%%%%%%%%%
\section{Numerical Examples}
\label{sec:simulations}
%%%%%%%%%%%%%%%%%%%%%%%%%%%%%%%%%%%%%%%%%%

\begin{figure}[thbp] % 'h' means "here" - adjust placement as needed (e.g., 't' for top, 'b' for bottom)
    \centering
    % Top row
    \begin{subfigure}[b]{0.23\textwidth} % Width of each subfigure (adjust as needed)
        \centering
        \includegraphics[width=\textwidth]{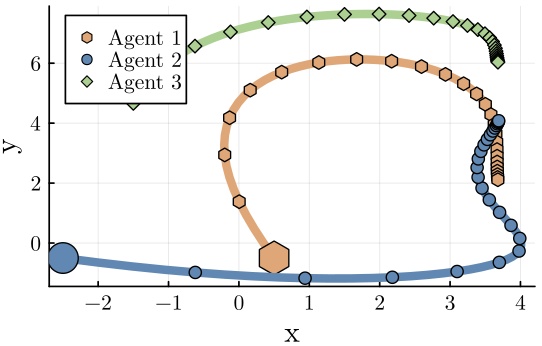} % Replace with your image file
        \caption{}
        \label{fig:consensus}
    \end{subfigure}
    \hfill % Adds horizontal space between subfigures
    \begin{subfigure}[b]{0.23\textwidth}
        \centering
        \includegraphics[width=\textwidth]{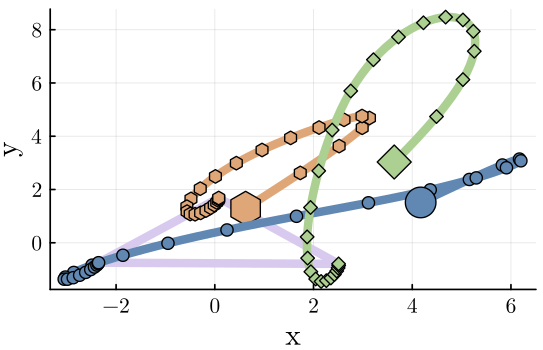}
        \caption{}
        \label{fig:formation}
    \end{subfigure}
    \hfill % Adds horizontal space between subfigures
    
    % Bottom row
    \begin{subfigure}[b]{0.23\textwidth}
        \centering
        \includegraphics[width=\textwidth]{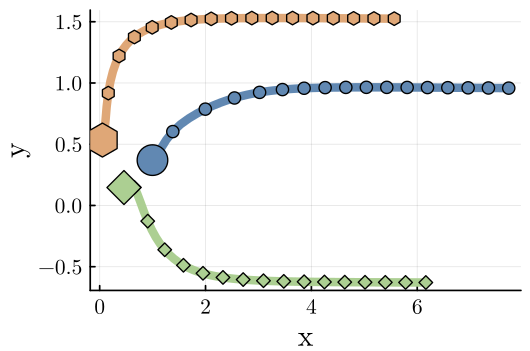}
        \caption{}
        \label{fig:flocking}
    \end{subfigure}
    \hfill
    \begin{subfigure}[b]{0.23\textwidth}
        \centering
        \includegraphics[width=\textwidth]{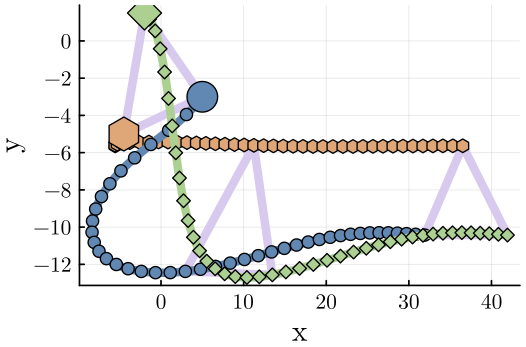}
        \caption{}
        \label{fig:moving}
    \end{subfigure}
    \caption{(a) Agents perform a hybrid goal of consensus in $x$ and tracking in $y$. (b) Agents perform the goal of reaching a triangular formation centered at the origin. (c) Agents perform the goal of flocking. (d) Agents perform a moving formation goal. In (c) and (d), the leader agent (1) tracks a constant rightward velocity vector.}
    \label{fig:experiments}
\end{figure}

To demonstrate the broad applicability of our framework, we implemented representative multi-agent optimal control examples of consensus, formation, and flocking.
We implemented our ADMM algorithm for nonlinear homological programs in the Julia programming language, using JuMP.jl with IPOPT to solve the individual agent subproblems and the conjugate gradient method from Krylov.jl to compute the projection operations \cite{JuMP, wachter2006implementation, Montoison2023krylov}. We implemented each example in a model predictive control setup, where at each timestep, the controllers for each agent solve the relevant optimal control homological program using Algorithm 1, implement the first control input, and repeat the process until convergence to the coordination goal.

The multi-agent system used for each example involved three vehicles obeying discrete double integrator dynamics in the plane ($d=2$), i.e. the Euler integration of \cref{eq:double-integrator}. Thus, each agent's state space is $\R^4$ and their control space is $\R^2$. Each dimension of each agent's control input was constrained to the range $[-2,2]$. Each agent's optimal control problem spanned a time horizon of $T=10$. We ran only 10 iterations of Algorithm 1 for each optimization step in MPC. This was sufficient for our MPC controllers to achieve convergence to the coordination goals while keeping the computation time low. In all experiments, each agent's initial positions and velocities were chosen randomly from the range $[-5,5]$. The code to reproduce the experiments can be found in our package AlgebraicOptimization.jl 

%---------------------
\paragraph{Consensus}
%---------------------
For this example, the agents are to reach consensus in the $x$-axis while reaching individual tracking goals in the $y$-axis. Thus each agent's subproblem is a linear quadratic tracking problem in $y$-space with a quadratic cost on the control activation. The coordination sheaf is defined over a fully connected communication topology, where each restriction map is simply the projection onto the $x$ coordinate. Edge potentials are the standard norm squared potential function encoding the consensus goal. The results of this controller run for 100 iterations are shown in Fig.~\ref{fig:experiments}.a.

%--------------------------------
\paragraph{Stationary Formation}
%--------------------------------
The goal in this example is for the agents reach a triangle formation centered at the origin. As such, each agent's subproblem is a standard linear quadratic regulation problem to drive the state to 0 with a quadratic penalty on control activation. The coordination sheaf is over a fully connected communication topology with the restriction maps projecting onto the position components of the state vector. The formation goal is encoded using edge potential functions of the form $U_e(y)=(1/2)\|y-b_e\|_2^2$ for desired $b_e$. The results of this controller run for 100 iterations are shown in Fig.~\ref{fig:experiments}.b.

%--------------------
\paragraph{Flocking}
%--------------------
For this example, agents implement the standard flocking goal of reaching consensus in velocities while staying a fixed distance away from all other agents. The constant sheaf $\underline{\R}^4$ on a fully connected communication topology along with potential functions summing the standard consensus potential function on the velocity components and the fixed distance potential function with $r^2=5$ on the position components. Each agents' objective is to minimize total control activation. Additionally, a designated leader agent tracks a constant rightward velocity vector. The results of this controller run for 65 iterations are shown in Fig.~\ref{fig:experiments}.c. Computing the distance between each agent confirms that they reached the desired pairwise distance of $\sqrt{5}$.
%-----------------------------
\paragraph{Moving Formation}
%-----------------------------
This example combines a formation goal in positions with a consensus goal in velocities. As such, the coordination sheaf is the constant sheaf $\underline{\R}^4$ on the three vertex path graph. This encodes a leader-follower topology with the middle agent in the path acting as the leader. The leader's objective is to track a constant rightward velocity vector and minimize its control actuation while the followers' objectives are to simply minimize control actuation. The edge potential functions are of the form $U_e(y)=(1/2)\|y-b_e\|_2^2$ where the velocity coordinates of each $b_e$ are 0 encoding consensus in velocity while the position coordinates are chosen to achieve a fixed displacement between the leader and its followers. The results of this controller run for 160 iterations are shown in Fig.~\ref{fig:experiments}.d.

%%%%%%%%%%%%%%%%%%%%%%%%%%%%%%%%%%%%%%%%%%
\section{Future Work}
\label{sec:conclusion}
%%%%%%%%%%%%%%%%%%%%%%%%%%%%%%%%%%%%%%%%%%

The results presented here can be extended in many directions. Having proposed a pipeline to translate coordination problems into homological programs and solve them with ADMM, a next step would be correct-by-construction software synthesis given user-specified coordination specifications.  Other promising areas for further developments include incorporating sheaf models of coordination into graph neural network architectures (see~\cite{zaghen_sheaf_2024, bodnar_neural_2022, hansen_sheaf_2020}) in order to solve agent-level or system-level tasks. While we attended primarily to convex potential functions, there is a rich theory to be developed for non-differentiable or non-convex potentials and the pursuit of local minima.

%%%%%%%%%%%%%%%%%%%%%%%%%%%%%%%%%%%%%%%%%%%%%%
\appendix
%%%%%%%%%%%%%%%%%%%%%%%%%%%%%%%%%%%%%%%%%%%%%%

\begin{proof}[Proof of \cref{thm:hp-convexity}]
    Because each $f_i$ is convex and the sum of convex functions is convex, the objective is clearly convex. It remains to show that the constraint set $\mathcal{C}=\{\mathbf{x}\mid L_\mathcal{F}^{\nabla U}\mathbf{x}=0\}$ is convex. For this, we utilize the fact that $L_\mathcal{F}^{\nabla U}$ is the gradient of the function $g(\mathbf{x}) = U(\delta_\mathcal{F} \mathbf{x})$. This function is clearly differentiable and is convex because the precomposition of a convex function with a linear map is convex (\cite{theRock}, Theorem 5.7). Thus we can apply Theorem 23.5 of \cite{theRock} to obtain the following
    \[
    \mathbf{x}'= \nabla g(\mathbf{x})\iff \mathbf{x}\in \partial g^*(\mathbf{x}'),
    \]
    where $\partial$ denotes the subdifferential and $g^*$ is the Legendre transform of $g$.
    Instantiating the above with $\mathbf{x}'=0$, we see that the set $\{\mathbf{x}\mid \nabla g(\mathbf{x})=0\}$ is equivalent to the set $\{\mathbf{x}\mid \mathbf{x}\in \partial g^*(0)\}$. The set of subgradients of a convex function at any point is closed and convex, therefore, the set $\mathcal{C}$ is closed and convex. So $\mathsf{P}$ asks to minimize a convex function restricted to a closed convex subset of its domain and is therefore a convex optimization problem.
\end{proof}

\begin{proof}[Proof of \cref{thm:nonlinear-projection}]
    For notational convenience, we let $L=L_\FF^{\nabla U}$ and $\delta=\delta_\FF$ for the remainder of the proof.
    First note that plainly $\mathbf{b}$ is the unique global minimizer of $U$, and $U(\mathbf{b})=0$. Thus $\nabla U(y)=0$ only when $y=\mathbf{b}$. We know that $L= \nabla(U\circ\delta)= \delta^\top\circ\nabla U \circ\delta$. This entails that $\dot{\mathbf{x}}$ is always in $\image\delta^\top$ and therefore orthogonal to $\ker\delta$. We thus decompose $\mathbf{x}$ as $\mathbf{x}^{||} + \mathbf{x}^\perp$ where $\mathbf{x}^{||}$ is the orthogonal projection onto $\ker\delta=H^0(G;\FF)$ and focus on the evolution of $\mathbf{x}^\perp$ under the dynamics restricted to $\image\delta^\top$.

    Let $\Psi(\mathbf{x}^\perp)=U(\delta \mathbf{x}^\perp)$. Because $\mathbf{b}$ is the unique minimizer of $U$ with $U(\mathbf{b})=0$, $\Psi$ vanishes only at $\delta \mathbf{x}^\perp =\mathbf{b}$. Moreover, $\delta \mathbf{x}^\perp = \mathbf{b}$ holds if and only if $\mathbf{x}^\perp = \delta^+\mathbf{b}$. To show the first direction of this implication, we have
    \[
    \delta \mathbf{x}^\perp=\mathbf{b} \implies \delta^+\delta \mathbf{x}^\perp = \delta^+\mathbf{b}\implies \mathbf{x}^\perp=\delta^+\mathbf{b},
    \]
    where the second implication holds because $\delta^+\delta$ is projection onto $\image \delta^\top$ and is therefore the identity on $\mathbf{x}^\perp$. For the reverse direction, we have
    \[
    \mathbf{x}^\perp =\delta^+\mathbf{b}\implies \delta \mathbf{x}^\perp=\delta\delta^+ \mathbf{b}\implies \delta \mathbf{x}^\perp = \mathbf{b},
    \]
    where the second implication holds because $\delta\delta^+$ is projection onto $\image \delta$ and is therefore the identity on $\mathbf{b}$.

    So, $\Psi(\mathbf{x}^\perp)$ vanishes only when $\mathbf{x}^\perp = \delta^+\mathbf{b}$ implying that $\nabla\Psi(\mathbf{x}^\perp)=L(\mathbf{x}^\perp)=0$ if and only if $\mathbf{x}^\perp=\delta^+\mathbf{b}$ and that $\Psi(\mathbf{x}^\perp)$ is globally positive definite about $\delta^+\mathbf{b}$. Furthermore, radial unboundedness of $U$ implies radial unboundedness of $\Psi(\mathbf{x}^\perp)$. Finally, we see that
    \[
    \dot{\Psi}=\langle\nabla\Psi(\mathbf{x}^\perp), \dot{\mathbf{x}}^\perp\rangle = \langle L(\mathbf{x}^\perp), -\alpha L(\mathbf{x}^\perp)\rangle \leq 0
    \]
    with equality only when $\mathbf{x}^\perp=\delta^+\mathbf{b}$. This shows the $\dot{\Psi}$ is globally negative definite around $\delta^+\mathbf{b}$, meaning $\Psi$ is a global Lyapunov function about $\delta^+\mathbf{b}$ for the dynamics in \cref{eq:nonlinear-diffusion-dynamics} restricted to $\image \delta^\top$. Therefore $\delta^+\mathbf{b}$ is globally asymptotically stable meaning $\mathbf{x}^\perp\to\delta^+ \mathbf{b}$. For the unrestricted dynamics in \cref{eq:nonlinear-diffusion-dynamics}, this means 
    \[\underset{t\to\infty}{\lim} \mathbf{x}(t) = \mathbf{x}^{||}(0) + \delta^+ \mathbf{b},\]
    which is orthogonal projection of $\mathbf{x}(0)$ onto $\delta^+\mathbf{b}+\ker\delta=\delta^+\mathbf{b} + H^0(G;\FF)$.

    To finish the proof, it remains to show that $\delta^+\mathbf{b} + H^0(G;\FF) = \ker L$. For the first direction, let $z\in \delta^+\mathbf{b} + \ker\delta$. Then we know $z$ is of the form $\delta^+\mathbf{b} + \mathbf{x}$ for some $\mathbf{x}$ such that $\delta \mathbf{x}=0$. We have
    \[
    \begin{array}{ll@{}ll}
    L(z) &= \delta^\top(\nabla U(\delta(\delta^+\mathbf{b}+\mathbf{x})))\\ &= \delta^\top(\nabla U(\delta\delta^+ \mathbf{b} +\delta \mathbf{x})) = \delta^\top(\nabla U(\mathbf{b}))=0,
    \end{array}
    \]
    where the third equality holds because $\delta\delta^+$ is projection onto $\image \delta$ and $\mathbf{x}\in\ker \delta$. For the reverse direction, let $L(\mathbf{x})=0$ for some $\mathbf{x}$. We know that this only holds when $\delta \mathbf{x}=\mathbf{b}$ which implies that $\delta^+\delta \mathbf{x} = \delta^+\mathbf{b}$. Because $\delta^+\delta$ is projection onto $\image\delta^\top = (\ker\delta)^\perp$, we know that $\mathbf{x}$ must be of the form $\delta^+\mathbf{\mathbf{b}} + \mathbf{x}'$ for some $\mathbf{x}'\in\ker \delta$. This completes the proof.
\end{proof}

\begin{proof}[Proof of \cref{thm:admm}]
    This setup satisfies the criterion for convergence in \cite[\S 3.2.1]{boyd_distributed_2010}.
\end{proof}

%%%%%%%%%%%%%%%%%%%%%%%%%%%%%%
\bibliographystyle{ieeetr}
\bibliography{tyler,hans}
%%%%%%%%%%%%%%%%%%%%%%%%%%%%%%

\end{document}